\documentclass[a4paper,12pt]{amsart} \usepackage{euscript,eufrak,verbatim}
\usepackage[final]{epsfig}
\def\w{\omega} \def\nn{\nonumber} \def\ms{m^*} \def\a{\alpha} \def\t{\theta}
\def\e{\varepsilon} \def\Xt{X_{\t}} \def\Gt{G_{\t}} \def\Bt{F_{\t}}
\def\Ut{U_{\t}} \def\Dt{D_{\t}} \def\Rt{R_{\t}}
\def\Jk{\mathcal{J}_k}
\def\U{\EuScript{U}}\def\sS{\mathcal{H}} 
\DeclareMathOperator{\diam}{diam} \DeclareMathOperator{\card}{Card}
\def\Z{\mathbb{Z}} \def\R{\mathbb{R}} \def\S{\mathbb{S}} \def\N{\mathbb{N}}
\def\Q{\mathbb{Q}} \def\Am{\EuScript{A}_{\mu,\w}} \def\Bm{\EuScript{B}_{\mu,t}}
\def\Ds{\EuScript{D}_s} \def\Es{\EuScript{E}_s} \def\Cs{\EuScript{C}_s}
\newtheorem{theorem}{Theorem}[section]
\newtheorem{lemma}[theorem]{Lemma}
\newtheorem{proposition}[theorem]{Proposition}
\newtheorem{corollary}[theorem]{Corollary}

\newtheorem{definition}[theorem]{Definition}
\DeclareMathSymbol{\varnothing}{\mathord}{AMSb}{"3F} \begin{document}

\title{Inhomogeneous Diophantine Approximation and Angular Recurrence for
Polygonal Billiards} \author{J\"org Schmeling} \author{Serge Troubetzkoy}
\address{Freie Universit\"at Berlin, FB Mathematik und Informatik, Arnimallee
2-6, D--14195 Berlin and Institut de math\'ematiques de Luminy, CNRS Luminy,
Case 907, F-13288 Marseille Cedex 9, France} \email{schmeling@math.fu-berlin.de}
\address{Centre de physique th\'eorique and Institut de math\'ematiques de
Luminy, CNRS Luminy, Case 907, F-13288 Marseille Cedex 9, France}
\email{troubetz@iml.univ-mrs.fr}
\urladdr{http://iml.univ-mrs.fr/{\lower.7ex\hbox{\~{}}}troubetz/} \date{}
 \begin{abstract} For a given rotation number we compute the
Hausdorff dimension of the set of well approximable numbers.  We use this result
and an inhomogeneous version of Jarnik's theorem to show strong recurrence
properties of the billiard flow in certain polygons.  \end{abstract} \maketitle

\pagestyle{myheadings}

\markboth{INHOMOGENEOUS APPROXIMATIONS AND BILLIARDS}{J\"ORG SCHMELING AND SERGE
TROUBETZKOY}

\section{Introduction}\label{sec1} In the past decade in four independent
articles it was observed that the billiard orbit of any point which begins
perpendicular to a side of a polygon and at a later instance hits some side
perpendicularly retraces its path infinitely often in both senses between the
two perpendicular collisions and thus is periodic.  The earliest of these
articles is a numerical work of Ruijgrok which conjectures that every triangle
has perpendicular periodic orbits~\cite{R}.  In 1992 Boshernitzan~\cite{B} and
independently Galperin, Stepin and Vorobets~\cite{GSV} proved that for any
rational polygon, for every side of the polygon, the billiard orbit which begins
perpendicular to that side is periodic for all but finitely many starting points
on the side.  Finally for an irrational right triangle Cipra, Hansen and Kolan
have considered points which are perpendicular to one of the legs of the
triangle.  They showed that for almost every such point the billiard orbit is
periodic \cite{CHK}.  Here the almost everywhere statement is with respect to
the length measure on the side considered.  The method of their proof in fact
implies a stronger result, namely that in any (generalized) parallelogram or
right triangle, for each direction $\t$, the set of points $\Bt$ whose orbit
starts in the direction $\t$ and never returns parallel to itself has Lebesgue
measure 0~\cite{GT}.  The result of Cipra {\it et.~al.}  on periodic orbits
follows from the above observation.  Their computer simulation indicated that
perhaps one can improve the almost everywhere and that is the starting point of
this research.  In Theorem~\ref{thm1} we show that in any (generalized)
parallelogram, for each direction $\t$, the set of points whose orbit starts in
the direction $\t$ and never returns parallel to itself has lower box dimension
at most one half.  A corollary to this theorem is that the set of points whose
orbit starts perpendicular to a leg of a right triangle whose orbit is not
periodic has lower box dimension at most one half.

We then turn to the question, whether there are directions $\t$ for which we can
improve the constant 1/2 in Theorem \ref{thm1}.  In Theorem~\ref{thm2} we show
that if $\mu > 1$ is the approximation order of $\t$ by $\a$, then in fact the
lower box dimension of $\Bt$ is at most $(\mu+1)^{-1}$.  We have two
applications of Theorem~\ref{thm2}.  We consider the set $\Cs$ of directions
$\t$ for which the set of points which never return parallel to themselves have
lower box dimension at most $s \in [0,1/2]$.  In Theorem~\ref{thmCs} we prove
that the Hausdorff dimension of $\Cs$ is at least $s/(1-s)$, the set $\Cs$ is
residual and has box dimension 1.  Fix a direction $\t_0$ and consider the set
$\Ds(\t_0)$ of angles $\a$ such that for any generalized parallelogram with
angle $\a$ the set of points in direction $\t_0$ which never return parallel to
themselves has dimension at most $s \in [0, 1/2]$.  We also consider the set
$\Es$ of right triangles for which the set of non-periodic points which are
perpendicular to a fixed leg of the triangle has dimension at most $s \in
[0,1/2]$.  In Theorem~\ref{thmDs} we prove that the Hausdorff dimension of $\Ds$
and $\Es$ are at least $2s$, these sets are residual and have box dimension 1.

The proofs of Theorems~\ref{thmCs} and \ref{thmDs} come from purely number
theoretic arguments.  The approximation of $t$ by $\w$ is a classical area of
research in number theory which is referred to as inhomogeneous Diophantine
approximation.  A classical result in this direction is the theorem of Minkowski
\cite{C} which states that if $t$ is not in the orbit of $\w$ then $\|t + p\w\|
< 1/(4p)$ has infinitely many integer solutions $p$ and the constant 1/4 can not
be improved in general.  Here $\| \cdot \|$ is the standard distance on $\S^1$.
In the spirit of Minkowski's theorem we consider the set $$\{t \in \S^1:
\|t+p\w\| < p^{-\mu} \text{ for infinitely many } p \in \N\}.$$ For any $\mu >
1$ in Theorem \ref{thmAm} we prove that the Hausdorff dimension of this set
is $\mu^{-1}$. Theorem \ref{thmCs} follows easily
from a similar statement for approximations along an arithmetic subsequence.

Another classical result in number theory is Jarnik's theorem on the Hausdorff
dimension of well approximable irrational numbers \cite{J}.  In the spirit of
Jarnik's theorem we consider the set 
$$\{\w \in \S^1:  \|t+p\w\| <
p^{-\mu}\text{ for infinitely many } p \in \N\}.$$ For any $\mu > 1$ Levesley
has proven that the Hausdorff dimension of this set is $2/(1+\mu)$ \cite{L}.  
The proof of Theorem \ref{thmDs} uses a similar inhomogeneous Jarnik result 
generalized to arithmetic subsequences.

The structure of the paper is as follows.  In Section~\ref{section2} we describe
the needed background results on dimension theory and billiards In
Section~\ref{section3} we state the purely number theoretic results and prove
the new results.  In Section~\ref{section4} we state and prove the new billiard
results.

Finally we remark that a generalization in a different direction of Cipra~{\it
et.~al.}'s results was obtained by one of us in~\cite{Tr}.


\section{Preliminaries}\label{section2}

\subsection{Dimension} Let $Y$ be a subset of $\R^n$.  Let $N(\epsilon)$ denote
the minimal number of $\epsilon$ balls needed to cover $Y$.

\begin{definition} For a subset $Y$ of $\R^n$, the lower box dimension of $Y$,
denoted by $\dim_{LB}Y$ is given by $$ \liminf_{\epsilon \to 0} \frac{\log
N(\epsilon)}{\log 1/\epsilon} \ .  $$ The upper box dimension $\dim_{UB}$ is
defined similarly, replacing the $\liminf$ by $\limsup$.  If $\dim_{UB}Y$ and
$\dim_{LB}Y$ both exist and are equal, we define the box dimension of $Y$ to be
this value, and write $\dim_{B}Y$ = $\dim_{UB}Y$ = $\dim_{LB}Y$.
\end{definition}

For a subset $U$ of $\R^n$, we let $\diam(U)$ denote the diameter of the set
$U$.

\begin{definition} Let $s \in [0,\infty]$.  The $s$-dimensional Hausdorff
measure ${\mathcal H}^s (Y)$ of a subset $Y$ of $\R^n$ is defined by the
following limit of covering sums:  $$ {\mathcal H}^s (Y) = \lim_{\epsilon \to
0}\left ( \inf \left\{\sum^\infty_{i=1} \left(\diam U_i\right)^s :  Y \subset
\bigcup^\infty_{i=1} U_i \text{ and } \sup_{i} \diam U_i \leq
\epsilon\right\}\right ).  $$ \end{definition}

It is easy to see that there exists a unique $s_0 = s_0 (Y)$ such that
\begin{equation} \label{eq:HD} {\mathcal H}^s (Y) = \begin{cases} \infty &
\mbox{ for } s < s_0\\ 0 & \mbox{ for } s > s_0\ .  \end{cases} \end{equation}

\begin{definition} The unique number $s_0$ given by Equation~\eqref{eq:HD} is
defined to be the Hausdorff dimension of $Y$ and is denoted by $\dim_H Y$.
\end{definition} Standard arguments give that for a subset $Y$ of $\R^n$, $$
\dim_H Y \leq \dim_{LB} Y\leq \dim_{UB} Y $$ There are examples which show that
these inequalities may be strict.

The box dimension can also be defined in terms of covering sums.  The only
change being that the covering sets all have equal diameter.  We note that in
order to estimate the box dimension, it suffices that the diameters of the
covering sets tend to $0$ along a geometric sequence.

Lastly, we define the Hausdorff dimension of a measure:  \begin{definition} Let
$\mu$ be a Borel probability measure on $X$.  Then the Hausdorff dimension of
the measure $\mu$ is defined by $$ \dim_H \mu = \inf_Y \left\{ \dim_H Y :  \mu
(Y) = 1\right\}.  $$ \end{definition}

We remark that the dimension of a measure is clearly always less than or equal
to the dimension of its (Borel) support.  There is a well known method of
computing lower bounds on the dimension of a measure or a set.  Let $\U(x,r)$
denote the ball of radius $r$ centered at $x$.

\begin{lemma}\label{lsy} If for some finite measure $\mu$ $$\liminf_{r \to 0}
\frac{\log \mu(\U(x,r))}{\log r} \ge s \text{ on a set of $\mu$--positive
measures } $$ then the dimension of the measure is at least $s$.  \end{lemma}

A survey of the methods and results in dimension theory can be found in
\cite{Falconer,Pesin}.


\subsection{Directional billiard transformation}\label{sec2}

Consider a polygon $Q \subset \R^2$.  A billiard ball, i.e.~a point mass, moves
inside $Q$ with unit speed along a straight line until it reaches the boundary
$\partial Q$, then instantaneously changes direction according to the mirror
law:  ``the angle of incidence is equal to the angle of reflection,'' and
continues along the new line.  For the moment we assume that $Q$ is convex.  We
will describe the billiard map in $Q$ as a transformation of the set $X$ of rays
which intersect $Q$.

We parameterize $X$ via two parameters.  The first parameter is the angle $\t$
between the ray and the positive $x$--axis.  To define the second parameter
consider the perpendicular cross section $\Xt$ to the set of rays whose angle is
$\t$.  The set $\Xt$ is simply an interval, thus the second parameter is then
the arc-length on $\Xt$ for a fixed orientation.  For a non-convex polygon we
must differentiate the portion of rays which enter and leave the polygon several
times.  The above construction can be done locally, yielding the set $\Xt$ which
consists of a finite union of intervals.  Let $w$ be the unnormalized length
measure on $\Xt$.

The billiard map $T:X \to X$ takes the ray, which contains a segment of a
billiard trajectory oriented by the direction of its motion, to the ray which
contains the next segment of this trajectory after the reflection in the
boundary.

A polygon $Q$ is called a generalized parallelogram if every side of $Q$ is
parallel to one of two fixed directions.  We suppose that one of these
directions is the direction of the positive $x$-axis.  
Then there is a unique $0
< \a < \pi/2$ such that the angle between any pair of sides is $0,\a$ or
$\pi - \a$.  We suppose throughout that $\a \not \in \mathbb{Q} \cdot \pi$,
that is that $Q$ is not a rational polygon.  We remark that the billiard in a
right triangle is equivalent to the billiard in a rhombus consisting of four
copies of the triangle via the process of unfolding.

\begin{figure}[t] 
\begin{center} 
\epsfxsize=170pt 
\epsffile{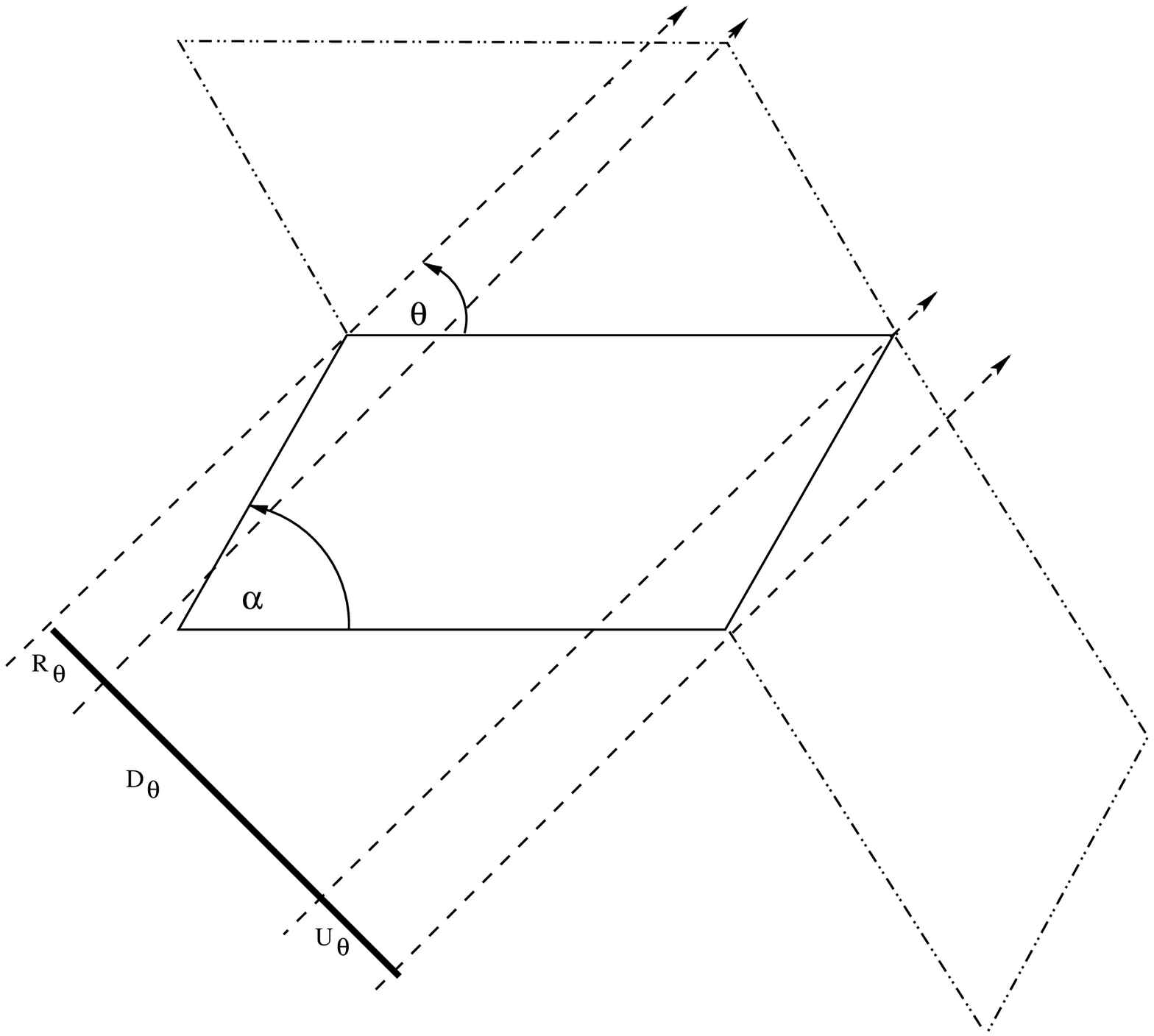} 
\caption{}\label{fig1} 
\end{center}
\end{figure}

The billiard map in a generalized parallelogram has the following special
properties (see figure~\ref{fig1}).  The interval $\Xt$ can be partitioned into
three sets $\Ut, \ \Rt$ and $\Dt$ which consist of a finite union of intervals
such that $T^2\Ut \subset X_{\t+2\a},\ T^2\Rt \subset \Xt$ and $T^2 \Dt \subset
X_{\t-2\a}$ and each interval in $\Ut, \ \Rt$ and $\Dt$ is mapped isometrically
onto its image.  The length measure $w$ is $T$--invariant.  Furthermore there is
a constant $K$ such that 
\begin{equation}\label{eq:small}
w(\Dt) \le K |\sin\t| \quad \hbox{ and } \quad w(\Ut) \le K |\sin(\t-\a)|.
\end{equation}
  We call the endpoints of the intervals of $\Ut$ and $\Dt$ the
singularities of $T^2$.

Let $Z_n := X_{(\t+2n\a)}$ and $Z := \cup_{n \in \Z} Z_n$.  The set $Z_n$ is
called the $n$th level of $Z$.  The set $Z$ is $T^2$-invariant with infinite
invariant measure $w$.

There are several good survey articles on billiards \cite{G1,G2,T}.  The
structure discussed here is described in more detail in \cite{GT}.


\section{Number theoretic results}\label{section3}

In this section we work with the circle $\R/\Z$ of unit length rather than the
circle of length $2\pi$.  The role of $\t$ and $\a$ in the billiard sections is
played by $t = (\t \mod \pi)/\pi$ and $\w = (2\a \mod \pi)/\pi$.

\subsection{Inhomogeneous Diophantine approximation} Let 
\begin{align*}
\Am & := \{t: \quad \|t+p\w\| < p^{-\mu} \text{ for infinitely many positive}\\
& \phantom{qqqqqqqqqqqqqqqqqqq}\text{ and infinitely
many negative } p \in \Z\}.
\end{align*}
Let $m,l \in \N$ be fixed such that $0 \le l < m$ and let
\begin{align*}
\Am(l,m) & :=\{t: \|t+p\w\| < p^{-\mu} \text{ for infinitely many positive
} \\
&\phantom{qqqqqqqq}\text{ and infinitely many negative } p\equiv l\mod m\}.
\end{align*}
Note that $\Am(l,m) \subset \Am$ and $\Am = \cup_{0 \le l < m}\Am(l,m)$.

Minkowski has shown that for $t$ not in the orbit of $\w$ the inequality
$\|t+p\w\| < 1/(4p)$ has infinitely many solutions and in general the constant
$1/4$ is optimal \cite{C}.  A simple Borel--Cantelli argument tells us that
directions which can be approximated better than in the statement 
of Minkowski's theorem have zero measure.

\begin{proposition}\label{borelcantelli} Suppose that $\mu > 1$ 
and $\w \not \in \Q$.  The Lebesgue measure of $\Am$ is 0 while 
$\dim_{B}\bigcap_{l=0}^{m-1}\Am(l,m) =1$.  The sets $\Am(l,m)$
are residual.  Hence, the set $\bigcap_{\mu \in \Q, \ \mu > 1}\bigcap_{l=0}^{m-1}
\Am(l,m)$ is also residual and has box
dimension~1.  \end{proposition}

\begin{proof} For simplicity, if $0 < a < b$ we denote the interval 
$(-b,-a)$ by $-(a,b)$.
Let $S(l,m)=\{n\in\Z\, :\, n\equiv l\mod m\}$, this yields
\begin{equation}\label{eAm} 
\Am(l,m) = \bigcap_{|k|=1}^{\infty}
\bigcup_{p \in S(l,m) \cap \rm{sgn} (k) [|k|m,\infty)} 
\left ( p\w-\frac1{2|p|^\mu},p\w+\frac1{2|p|^\mu}\right )
\end{equation} 
and 
$$\Am = \bigcap_{|k|=1}^{\infty}
\bigcup_{p \in \rm{sgn}(k) [|k|,\infty)}
\left ( p\w-\frac1{2|p|^\mu},p\w+\frac1{2|p|^\mu}\right )$$
thus since $\{p\w:  \ p\in S(l,m) \cap [|k|m,\infty)\}$ and $\{p\w:  \ p\in S(l,m) \cap -[|k|m,\infty)\}$ are dense 
$\Am(l,m)$ is residual. Thus the set
$\bigcap_{k=0}^{m-1}\Am(l,m)$ is dense, hence it's lower and upper box 
dimension are 1.  To see that the Lebesgue measure of $\Am$ is 0 note
that the sequence
$(2|p|^{-\mu})_{p \in \Z}$ has finite sum for $\mu > 1$.  Hence, by the
Borel--Cantelli lemma almost no point $t$ is contained in more than a finite
number of intervals $(p\w-1/2|p|^{\mu},p\w+1/2|p|^{\mu})$.  \end{proof}

Due to this proposition we consider the dimension of the set $\Am$ of $t$'s with
better approximations.

\begin{theorem}\label{thmAm} For any $\mu > 1$ and $m\in\N$ $$\dim_{H} \bigcap_{l=0}^{m-1}\Am(l,m) = \frac1{\mu}$$
\end{theorem}

\begin{proof} We begin by proving the upper bound.  For any $\e > 0$
Equation~\eqref{eAm} implies 
$$\sS^{\frac1{\mu}+\e }(\Am(l,m)) \le
2 \times \liminf_{k\to\infty}\sum_{p=k}^\infty \left ( \frac 1{p^\mu} \right
)^{\frac1{\mu}+\e} =0.  $$ We turn to the lower bound.  
For $n_k\in\Z$  let $\Jk := \rm{sgn}(n_k)[|n_k|,2|n_k|]$ 
and
$$\hat{T}_{n_k}^\mu:=\bigcup_{qm+k\in \Jk} \left
((qm+k)\w-\frac1{2(2|n_k|)^\mu},(qm+k)\w+\frac1{2(2|n_k|)^\mu}\right ).$$ Given any integer sequence
$\{n_k\}$ which meets all residue classes infinitely often for both positive and negative $n_k$ we have that $$\hat{S}_{\{n_k\}} := \bigcap_{k=1}^\infty
\hat{T}_{n_k}^\mu \subset \bigcap_{l=0}^{m-1}\Am(l,m).$$ To prove the lower bound we will use induction
to construct a sequence $n_k$ which is simply a subsequence of the denominators
$q_m$ of the continued fraction approximants.  Given the sequence $\{n_k\}$ we
will construct a measure $m$ supported on the set $\hat{S}_{\{n_k\}}$ whose
dimension is $\mu^{-1}$.  The idea of our construction is as follows:  the
faster the sequence $n_k$ increases the better the points 
$((qm+k)\w)_{qm+k\in \Jk}$ are distributed along the 
circle.  
If a point $x$ belongs to $\Am$ then
it has to be infinitely often in intervals of the form
$J_q=(q\w-|q|^{-\mu}/2,q\w+|q|^{-\mu}/2)$.  So for some large ${n_k}$ it has to fall
into $J_q$ for some $q \in \rm{sgn}(n_k) [0,2|n_k|$.  The set $S_{\{n_k\}}$ forces the point
to be in the ``right half'' of the collection of intervals $J_q$, where $q$ runs
from 0 to $2{n_k}$.  If the sequence $\{n_k\}$ is sufficiently sparse this
restriction turns out to be mild enough to maintain the dimension of $\Am$.

Let $n_1 = q_1$ and assume that the sequence $n_i$ is constructed up to $i =
k-1$.  Then the set $\bigcap_{l=1}^{k-1} \hat{T}_{n_l}^\mu$ consists of a finite
number of open intervals $\{\hat{I}^{k-1}_i\}$ satisfying \begin{equation}
\label{length} \left | \hat{I}^{k-1}_l \right | \le (2{|n_{k-1}|})^{-\mu}
\end{equation} with equality if and only if they are completely contained in
$\cap_{l=1}^{k-2} \hat{T}_{n_l}^\mu$.  We will drop the hat from the notation
for any interval for which equality holds.

We consider the set $$T^{\mu}_{n_k} := \bigcup_l I^k_l$$ and $$S_{\{n_k\}} :=
\bigcap_{k=1}^{\infty} T^{\mu}_{n_k}.$$ Clearly $S_{\{n_k\}} \subset
\hat{S}_{\{n_k\}}$.

Assume that inductively $\{n_i\}$ is already constructed for $i \le k-1$.
It is well known that for any $k\in\N$ the sequence 
$\{(pm+k)\w:  \ p \in \Z\}$ is well distributed in
$\S^1$, i.e.~for any continuous function $f$ and any $\e > 0$ there is a 
number $N$ such for any 
$q \in \Z$ $$\left | \frac1n \sum_{p=0}^{n-1} f((q+pm+k)\w) - \int
f \right | < \e$$ for any $n \ge N$.  
Thus if $k \mod 2m \le m-1$ we can choose $n_k$ to be the 
least positive $q_r > |n_{k-1}|$ such that the following two conditions hold:
\begin{equation}\label{well} \frac{|I^{k-1}_l|}{2} \le \frac{\card \{pm+k \in
\Jk :  \ (pm+k)\w \in I^{k-1}_l\}}{|n_k|} \le 2|I^{k-1}_l| \end{equation}
and \begin{equation}\label{verywell} \frac{\log\prod_{i=1}^{k-1}{|n_i|}}{\log
{|n_k|}}\to 0.  \end{equation} otherwise if $k \mod 2m \ge m$ 
we choose $n_k=-q_r$ satisfying Conditions 
\eqref{well} and \eqref{verywell}.  
This finishes the construction of the 
sequence
$\{n_k\}$. Clearly, the sequence $\{n_k\}$ meets all residue classes infinitely often for positive and negative $n_k$.

We next construct a measure on the set $S_{\{n_k\}}$ with the desired 
dimension.
We assume that $I^0$ consists of the whole circle and thus has length 1.  We
begin by a recursive definition of an outer measure $\ms$.  We put $\ms(I^0) :=
1$ and $$\ms(I^k_l) := \frac{\ms(I^{k-1}_j)} {\card\{j:  \ I^k_l \subset
I^{k-1}_j\}}.$$

We will now compute an upper bound on the outer measure.  Using
Equations~\eqref{well} and \eqref{length} we have \begin{eqnarray}\label{upper}
\ms(I^{k+1}_j) & \le & \frac{2\ms(I^k_l)}{|I^k_l|{|n_{k+1}|}}\nn\\ & = & 2
\ms(I^k_l) \frac{(2{|n_k|})^\mu}{{|n_{k+1}|}}\nn\\ & \le & 2^{(1 + \mu)k}
\prod_{i=1}^k \frac{({|n_i|})^\mu}{{|n_{i+1}|}}\\ & \le & 2^{(1 + \mu)k}
\frac{(n_1)^\mu}{{|n_{k+1}|}} \prod_{i=2}^k ({|n_i<})^{\mu-1}.\nn \end{eqnarray} The
fast growth rate on $n_k$ assumed by Equation~\eqref{verywell} implies that
$m(I^k_j) \to 0$ as $k \to \infty$.  Thus it is clear that $\ms$ satisfies
Kolmogorov's compatibility conditions and can be extended to a measure $m$ on
$S_{\{n_k\}}$.

We are now ready to prove that the dimension of the measure is at least $1/\mu$.
Equation~\eqref{upper} implies \begin{eqnarray*} \frac{\log m(I^{k+1}_l)}{\log
|I^{k+1}_l|} & \ge & \frac{\log 2^{(1+\mu)k} + \sum_{i=2}^k \log {|n_i|}^{\mu-1} +
\mu \log n_1 - \log {|n_{k+1}|}}{\log(2{|n_{k+1}|})^{-\mu}}\\ & = & \frac{-kC_1 -
C_2 \sum_{i=2}^k \log {|n_i|} -C_3 + \log {|n_{k+1}|}} {C_4 + \mu \log {|n_{k+1}|}}
\end{eqnarray*} where the $C_i$ are all positive constants.
Equation~\eqref{verywell} implies that $$ \lim_{k \to \infty} \frac{\log
m(I^{k+1}_l)}{\log |I^{k+1}_l|} \ge \frac1{\mu}.  $$ To use Lemma~\ref{lsy} to
conclude our result we must evaluate the ratio $\log m(I)/ \log |I|$ for the
``intermediate'' intervals $I$, that is when the interval $I$ satisfies
$$\frac1{(2|n_k|)^\mu} < |I| < \frac1{(2|n_{k-1}|)^\mu}.$$ Let $\mathcal{U}_{a}(B)$
be the neighborhood of radius $a$ around the set $B$ and let 
$r(I) := \card \{qm+k : \ qm+k \in \Jk \text{ and } (qm+k)\omega \in I \backslash
\mathcal{U}_{(2n_{k})^{-\mu}}(\partial I)\}.$

If $r=1$ then we denote by $I^k_l$ the unique such interval which is contained
in $I$.  Then $$\frac{\log m(I)}{\log |I|} = \frac{\log m(I^k_l)}{\log |I|} >
\frac{\log m(I^k_l)}{\log|I^k_l|}.$$

Now suppose $r \ge 2$.  Since $|n_k|$ is some continued fraction convergent
$q_{m_k}$ we have \begin{equation}\label{eee} \min_{|n_k| \le |p_1| < |p_2| \le |2n_k|}
d(p_1\omega,p_2\omega) \ge \frac1{|n_k|+2}.  \end{equation} We remark that this is the
only place where we use that the $n_k$ is a subsequence of the continued
fraction convergents.  Equations~\eqref{eee} and \eqref{length} imply (assuming
that $|n_k| \ge 2$) that $$|I| \ge (r-1) \frac1{|n_k|+2} \ge (r-1)
|I^k_l|^{\frac1{\mu}}.$$ Thus we can make the following estimate
\begin{eqnarray}\label{inter} \frac{\log m(I)}{\log |I|} & \ge & \frac{\log
m(I)}{\log ((r-1) \cdot |I^k_l|^{\frac1{\mu}})}\nonumber\\ & \ge & \frac{\log (r
\cdot m(I^n_l))}{\log ((r-1) \cdot |I^k_l|^{\frac1{\mu}})}\\ & = & \frac{\log r
+ \log m(I^k_l)}{\log (r-1) + \frac1{\mu}\log|I^k_l| } \nonumber \end{eqnarray}
For any $\e > 0$, using Inequality~\eqref{inter} we can choose $k(\e)$
sufficiently large such that for any $k \ge k(\e)$ we have that $$\frac{\log
m(I)}{\log |I|} \ge \frac{\log r + \frac{1}{\mu + \e}\log |I^k_l|}{\log(r-1) +
\frac1{\mu} \log |I^k_l|}.$$ Thus $$\frac{\log m(I)}{\log |I|} \ge
\frac{\mu}{\mu+\e} > \frac{1}{\mu}$$ where the last inequality holds for any
sufficiently small $\e$.

Now Lemma~\ref{lsy} implies our result.  \end{proof}


\subsection{Jarnik's theorem} 

Let 
\begin{align*}
\Bm & := \{\w:  \|t+p\w\| <
p^{-\mu}\text{ for infinitely many positive}\\
& \phantom{qqqqqqqqqqqqqqqqqqq}\text{ and infinitely
many negative } p \in \Z\}.
\end{align*}
Let $m,l \in \N$ be fixed such that $0\le l<m$ and let
\begin{align*}
\Bm(l,m) & := \{\w:  \|t+p\w\| < p^{-\mu} \text{ for infinitely many positive}\\
&  \phantom{qqqqqqqqqqqqqq}\text{ and infinitely many negative } 
p \equiv l\mod m\}.
\end{align*}
Then $\Bm(l,m)\subset\Bm$ and $\Bm\supset\bigcup_{0\le l<m}\Bm(l,m)$.

\begin{proposition}\label{borelcantelli1} Suppose that $\mu > 1$, and 
$\w \not \in \Q$.  
The Lebesgue measure of $\Bm(l,m)$ is 0 while $\dim_{B}\Bm(l,m) =1$.
The sets $\Bm(l,m)$ are residual.  Hence, the set 
$\bigcap_{\mu \in \Q, \ \mu>1}\bigcap_{0\le l <m}\Bm(l,m)$ is also
residual and has box dimension~1.  \end{proposition}

\begin{proof}

Let $S(l,m)=\{n\in\N\, :\, n\equiv l\mod m\}$, this yields
\begin{equation}\label{eBm1} \Bm(l,m) = \bigcap_{|k|=1}^{\infty} \bigcup_{p\in
S(l,m)\cap \rm{sgn}(k)[|k|,\infty)} \bigcup_{i=0}^{p-1} \left ( \frac{t+i}{p}
-\frac1{2|p|^{\mu+1}}, \frac{t+i}{p}+\frac1{2|p|^{\mu+1}}\right ) 
\end{equation}
hence $\Bm(l,m)$ is residual.  The set $\Bm(l,m)$ is dense, thus it's lower and
upper box dimension are 1.  We have $\sum 1/2p^{\mu+1} < \infty$, hence, by the
Borel--Cantelli lemma almost no point $\w$ is contained in more than a finite
number of intervals $(\frac{t+i}{p} -\frac1{2|p|^{\mu+1}},
\frac{t+i}{p}+\frac1{2|p|^{\mu+1}})$, i.e.~the Lebesgue measure of $\Bm(l,m)$ 
is 0.

\end{proof}

The following generalization of Jarnik's classical theorem \cite{J} is a 
special case of a more general result which has been proven by Levesley:  
\begin{theorem}\label{thmBm} \cite{L} For any $\mu > 1$
$$\dim_H \Bm = \frac{2}{1+\mu}.$$ \end{theorem}

We are going to improve Levesley's Theorem in this special setting.  The
proof essentially follows his ideas and we do not claim any originality.
\begin{theorem}\label{general}  For any $\mu > 1$ and $m \in \N$ $$\dim_H
\bigcap_{l=0}^{m-1}\Bm(l,m) = \frac{2}{1+\mu}.$$ \end{theorem}

We will need the (adapted to our case) notion of an {\it ubiquitous} system.
For fixed $l,m$ let $R_k(l,m)=\{t \in [0,1)\, :\, \|(km+l)t-\w\|=0\}$.

\begin{definition} Let $\rho\colon\N\to\R$ with $\lim_{N\to\infty}\rho(N)=0$ be
given.  Then the family of point sets $\{ R_k(l,m) \}$ is said to be {\em
ubiquitous} with respect to $\rho$ if \[
\lim_{N\to\infty}Leb\left([0,1)\setminus\bigcup_{k=1}^N\U(R_k(l,m),\rho(N))\right)=0.
\] \end{definition} We can express the sets $\Bm(l,m)$ in the following way \[
\Bm^+(l,m):=\bigcap_{K=0}^\infty\bigcup_{k=K}^\infty \U\left(R_k(l,m),
(km+l)^{-\mu-1}\right),
\]
\[
\Bm^-(l,m):=\bigcap_{K=0}^\infty\bigcup_{-k=K}^\infty \U\left(R_k(l,m),
(km+l)^{-\mu-1}\right)  \]
and
\[
\Bm(l,m)=\Bm^+(l,m)\cap\Bm^-(l,m).
\]
 We are going to use the following special case of
the more general Theorem 2 in~\cite{DRV} proved by Dodson, Rynne and Vickers
\begin{theorem}[\cite{DRV}]\label{dpv} Suppose that for each $0\le l<m$ the
family $\{R_k(l,m)\}$ is ubiquitous with respect to $\rho$.  Then \[
\dim_H\bigcap_{l=0}^{m-1}\Bm(l,m)\ge
\min\left\{1,\limsup_{N\to\infty}\frac{\log\rho(N)}{\mu+1}\right\}.  \]
\end{theorem} 

Levesley proved by using discrepancy estimates that the
system $\{\bigcup_{l=0}^{m-1}R_k(l,m)\}$ is ubiquitous with respect to 
the function
$\rho(N):=K\frac{\log^{5+\epsilon}N}{N^2}$.  Noting that the sequence $(k m t
+l)_k$ has the same discrepancy as the sequence $(k m t)_k$ we can follow the
arguments of Levesley and obtain that $\{R_k(l,m)\}$ is ubiquitous with respect
to the function
$\rho(N):=K\frac{\log^{5+2\e}(Nm+l)}{(Nm+l)^2}$.  Now Theorem~\ref{dpv} implies
Theorem~\ref{general}.

\section{Billiard results}\label{section4} \subsection{Strong recurrence in all
directions}

Fix a generalized parallelogram.  For every direction $\t$ let $\Bt$ be the set
of $x \in \Xt$ whose forward orbit never returns parallel to $x$ and let $\Gt :=
\Xt \setminus \Bt$.

\begin{theorem}\label{thm1} For any generalized parallelogram, for every
direction $\t$, the set $\Bt$ has lower box counting dimension at most 1/2.
\end{theorem}

This is an improvement of the main theorem of \cite{GT} which asserts that 
for every $\t$ the set $\Bt$ has measure 0. We remark that the set $\Gt$ 
is open and dense for each $\t$. 

Unfolding a right triangle with angle $\a/2$ around its right angle yields a 
rhombus with angle $\a$.  We consider
a direction $\t$ in the rhombus which is perpendicular to one of the legs of the
triangle.  The orbits which start in $\Xt$ and return to $\Xt$, when considered
as orbits in the triangle are twice perpendicular to a side.  As explained in
the introduction, such orbits must be periodic.  Thus applying
Theorem~\ref{thm1} to this direction yields the following improvement of the
theorem of Cipra, Hansen and Kolan \cite{CHK}:

\begin{corollary} For a right triangle, the set of points which are
perpendicular to one of the legs of the triangle and whose orbit is not periodic
has lower box counting dimension at most 1/2.  \end{corollary}

To prove Theorem~\ref{thm1} we start with a lemma which is essentially contained
in \cite{GT}

\begin{lemma}\label{lemma1} Fix an interval $\Xt$ and a positive integer $N$.
Then there is a partition of
$\Ut$ into $j_N \le CN$ intervals for some positive constant $C$ such the
forward orbit of each interval of the partition is an isometric mapping until
one of the following happens 
1) the orbit of the interval reaches level $N$
before returning to $\Xt$, or 
2) the orbit of the interval returns to $\Xt$ without having reached
level $N$.  The orbits of all these intervals are mutually disjoint until they
possibly return to $\Xt$.  \end{lemma}

A similar statement is true for $\Dt$ and level $-M$.

\begin{proof} Consider the set $V_N$ of the singularities of $T^2$ on levels 0,1
to $N$.  The set $V_N$ has cardinality at most $CN$ (here $C$ can be taken to be
twice the number of vertices of the polygon).  For each $v \in V_N$ consider the
first (if any) preimage which is on level 0.  We remark the length of time for
$v \in V_N$ to be mapped to level 0 is not necessarily bounded.  This forms the
indicated partition of $\Ut$.  The mapping $T^{2j}$ is continuous on each
element of the partition until it is mapped into level $N+1$ or returns to 
level 0.  

The disjointness follows from the invertability of $T^2$.
To see this suppose that $x,y \in \Xt$, $ x \ne y$ but $T^{2j}x =
T^{2k}y$ with $j \ge k$. We can not have $j=k$ since then the point $z$ has
two preimages. Since $T^{2(j-k)}x = y$ the orbit of $x$ has returned 
to $\Xt$ at time $2(j-k) > 0$ and
its orbit is disjoint from $y$'s before this time as claimed.

If one of the intervals (call it $J$) never returned to $\Xt$ or reached level 
$N$ then the set $\cup_{j \in \N} T^{2j}J$ would consists of an
infinite number of disjoint intervals of the same length, yet it
would be contained in $\cup_{0 \le n \le N} Z_n$.
This later set consists of a finite union of intervals of finite length
yielding the desired contradiction.
\end{proof}

\begin{proof}[Proof of Theorem~\ref{thm1}] Fix an arbitrary $\t$.  
We will argue about the dimension of $\Bt \cap \Dt$.
The argument for $\Bt \cap \Ut$ is
similar with the proviso that since we use Equation \eqref{eq:small}
$\t - 2N\a$ must be close to 0 for $\Dt$ while for $\Ut$, $\t + 2N\a$
must be close to $\a$. 
Because of this proviso, since we need both $\Bt \cap \Dt$ and
$\Bt \cap \Ut$ to be simultaneously small,
rather than using the approximants by $2\a$
we will use the approximants by $\a$, using the even iterates for
the set $\Bt \cap \Dt$ and the odd iterates for the set $\Bt \cap \Ut$.  

Consider the set $F_N$ of those points $x \in \Dt$ whose forward orbit reaches
$X_{\t-2(N+1)\a}$ before returning to $\Xt$.  We call this event $x$ 
exiting from
level $-N$ (necessarily to level $-(N+1)$).  When the orbit of $x$ reaches 
level $-N$
and is in the process of exiting level $-N$ it must pass through the ``gate''
$D_{\t-2N\a}$ on level $-N$.

Using Lemma~\ref{lemma1} clearly at most $CN$ of the intervals 
exit before returning to
$\Xt$.  Hence the set $F_N$ consists of at most $CN$ intervals which we call
$I_i$.  When they exit level $-N$ before returning to level 0 the sum of there
lengths can not be more than the width of the gate on level $-N$ through which
they exit.  Here we emphasize that we use the statement in Lemma~\ref{lemma1}
that the exiting orbits are all disjoint.  Thus the total length of the
intervals exiting level $-N$ is at most the total length of the gate of level
$-N$, i.e.~at most $K|\sin(\t - 2N\a)|$.

We use $[x]$ to denote the integer part of $x$.  Let $\{a_i\}$ be a sequence of
positive numbers and let $n \in \N$.  A simple estimate yields (see \cite{KS}):
\begin{equation}\label{e1} \sum_{k=1}^n
\left(\left[\frac{na_k}{\sum_{i=1}^na_i}\right] + 1 \right) \leq 3n \ .
\end{equation}

The $1/2+\epsilon$ covering sum of $\Bt \cap \Dt$ can be estimated by covering
$F_N$ by intervals of the ``average exiting length'' i.e.~by intervals of length
$K|\sin(\t - 2N\a)|/j_N$.

The total number of intervals of average length which we need to cover is at
most $$\sum_{k=1}^{j_N}\left(\left[\frac{|I_k|(j_N)} {\sum_{i=1}^{j_N}|I_i|}
\right] +1 \right).$$

Using Inequality~\eqref{e1} and the fact that $j_N \le CN$, we have
\begin{align}\label{eq:ineq} \sS^{1/2+\e} (\Bt \cap \Dt) & \le
\sum_{k=1}^{j_N}\left(\left[\frac{|I_k|(j_N)}{\sum_{i=1}^{j_N}|I_i|} \right] +1
\right) \left(\frac{K(|\sin(\t - 2N\a)|)}{j_N}\right)^{1/2+\e} \nonumber\\
& \leq
\frac{3j_N}{(j_N)^{1/2+\e}} \left (K(|\sin(\t - 2N\a)| \right )^{1/2+\e} \\ & \le
LN^{1/2-\e}\left((\t-2N\a)\mod \pi \right)^{1/2+\e}\nonumber \end{align} for some
constant $L$.  To get an effective cover we just need to find a sequence 
$(N_k)$
such that the inhomogeneous Diophantine approximation $\|(\t -2N_k\a)/\pi\|$ is
small.  It is well known that for any $\delta>0$ the inequality
\begin{equation}\label{deltadio}  \|(\t - 2p\a)/\pi\|<p^{-1+\delta} \end{equation}
is solvable for infinitely many $p\in\N$.  We remark that this is a much weaker
statement than Minkowski's Theorem \cite{C}.  Hence, we can find sequences
$(N_k)_{k \in \N}$ such that \begin{align}\label{s1/2} \sS^{1/2+\e}(\Bt \cap
\Dt)&\le LN_k^{1/2-\e}\left(N_k^{-1+\delta} \right)^{1/2+\e}\nonumber\\&\le
LN_k^{-2\e+(\delta/2)+\delta\e}.  \end{align} For $1>\e>\delta>0$ the
right--hand--side tends to zero as $k\to\infty$.

By the observation that we covered $\Bt \cap \Dt$ by equal length intervals, 
we have
$\dim_{LB}(\Bt \cap \Dt) \leq 1/2$.  

To estimate the dimension  of $\Bt \cap \Ut$
we need to use the second part of Equation \eqref{eq:small} in  
Equation \eqref{eq:ineq}. With this change it is useful to consider 
the variable
$\t' := \t - \a$.  Noticing that $\t + 2p\a - \a = \t' + 2p\a$ and
remarking that $\|((\t' + 2p\a)/\pi\| < p^{-1 + \delta}$ has infinitely 
many integer solutions we conclude that $\dim_{LB}(\Bt \cap \Ut) \leq 1/2$
and thus that  $\dim_{LB}(\Bt) \leq 1/2$.
\end{proof}


\subsection{Number theory implies stronger recurrence in some directions}

In the previous section we proved a recurrence theorem for all directions $\t$.
For our first results in this section we consider special directions $\t$ which
are extremely well approximable.  For these directions we can conclude a
stronger recurrence result.

\begin{theorem}\label{thm2} Consider a generalized parallelogram with angle 
$\a$ and let $\t$ be a direction such that 
$\|(\t-p\a)/\pi\| < p^{-\mu}$ for infinitely many $p \in 2\N$
and
$\|(\t + p\a)/\pi\| < p^{-\mu}$ for infinitely many $p \in 2\N+1$.
Then the set $\Bt$ has lower box counting dimension at most $1/(\mu +1)$.
\end{theorem}

\begin{proof} The proof is essentially identical to the proof of
Theorem~\ref{thm1}.  One only needs to change Inequality~\eqref{deltadio} to the
assumption of the theorem and change~\eqref{s1/2} accordingly.  \end{proof}

We apply the above result to perpendicular orbits in right triangles.  Here we
assume that the base of the triangle is parallel to the $x$--axis.

\begin{corollary} Fix $\mu > 1$ and consider a right triangle with angle 
$\a/2$, such that the approximations of $\pi/2$ by the orbit of $\a$ in
the sense of the previous theorem is of order $\mu$.
Then the set of points which are perpendicular to the base of the triangle and
whose orbit is not periodic has lower box counting dimension at most $1/(\mu
+1)$.  \end{corollary}

Theorem~\ref{thm2} indicates that to get a better result than Theorem~\ref{thm1}
one needs to consider directions with better approximations, i.e.~$\mu > 1$.
Let $$\Cs := \{\t:  \dim_{LB}(\Bt) \le s\}.$$

\begin{theorem}\label{thmCs} Fix an arbitrary generalized parallelogram.  For
any $s \in [0,1/2]$ we have $$\dim_{H}\Cs \ge \frac{s}{1-s},$$ $\Cs$ is residual
and has box dimension 1.  \end{theorem}

\begin{proof} The first statement is an immediate consequence of the previous
theorem and Theorem~\ref{thmAm} with $m=2$.  
The second statement is an immediate
corollary of Propositions ~\ref{borelcantelli}.  \end{proof}

Next we will fix the direction and vary the generalized polygons.  To do this we
assume that one of the sides of the polygon is fixed and we measure the
direction $\t$ with respect to this fixed direction.  We remark that our
estimates on $\Bt$ will depend only on the angle $\a$ and no other parameters of
the generalized parallelogram.  Fix a direction $\t_0$.  Consider the sets
\begin{eqnarray*} \Ds & := & \{\a \in \S^1:  \dim_{LB}(F_{\t_0})<s
\text{ for all generalized parallelograms}\\
&&\phantom{1111111111111111111111111111111111} \text{ with angle } \a \} \quad
\text{and} \\ \Es & := & \{\a \in \S^1:  \dim_{LB}(F_{perp}) <s \text{ for the
right triangle} \\
&&\phantom{1111111111111111111111111111111111} \text{ with angle } \a/2\} \end{eqnarray*} where perp is a direction
perpendicular to one of the legs of the triangle.

\begin{theorem}\label{thmDs} (a) For  $s \in [0,1/2]$ we have
$\dim_H\Es \ge 2s.$ The set $\Es$ is residual and has box
dimension 1. \\
(b) For $\t_0$ fixed and $s \in [0,1/2]$ we have
$\dim_H\Ds \ge 2s$. 
The set $\Ds$ is residual and has box dimension 1. 
 \end{theorem}

\begin{proof} The results on Hausdorff dimension
are immediate consequences of Theorem~\ref{thm2} 
and Theorem~\ref{general} with $m=2$.  
The residuality and box dimension results
are immediate corollaries of
Proposition~\ref{borelcantelli1}.   \end{proof}

\bigskip

\noindent {\bf Acknowledgements} We thank Martin Schmoll and Anatole Stepin for
useful suggestions.


\begin{thebibliography}{99} \footnotesize{ \bibitem[B]{B} M.~Boshernitzan, {\it
Billiards and rational periodic directions in polygons,} Amer.~Math.~Monthly 99
(1992) 522--529.

\bibitem[C]{C} J.~W.~S.~Cassels, \emph{An introduction to Diophantine
approximation}, Cambridge University Press, Cambridge, 1957.


\bibitem[CHK]{CHK} B.~Cipra, R.~Hanson, and A.~Kolan, {\it Periodic trajectories
in right triangle billiards}, Phys.~Rev.~E 52 (1995) 2066-2071.

\bibitem[DRV]{DRV} M.  Dodson, B.  Rynne, and J.  Vickers, {\it Diophantine
approximation and a lower bound for Hausdorff dimension}, Mathematika 37 (1990),
59--73


\bibitem[F]{Falconer} K.~Falconer.  \emph{The geometry of fractal sets,}
Cambridge University Press, Cambridge, 1985.

\bibitem[GSV]{GSV} G.~Galperin, A.~Stepin and Ya.~Vorobets, {\it Periodic
billiard trajectories in polygons:  generation mechanisms,} Russian
Math.~Surveys 47 (1992) 5--80

\bibitem[G1]{G1} E.~Gutkin, {\it Billiard in polygons,} Physica D, 19 (1986)
311-333.

\bibitem[G2]{G2} E.~Gutkin, {\it Billiard in polygons:  survey of recent
results,} J.~Stat.~Phys., 83 (1996) 7--26.

\bibitem[GT]{GT} E.~Gutkin, and S.~Troubetzkoy, {\it Directional flows and
strong recurrence for polygonal billiards,} in Proceedings of the International
Congress of Dynamical Systems, Montevideo, Uruguay, F.~Ledrappier et.~al.~eds
(1996) 21--45.

\bibitem[J]{J} V.~Jarnik, {\it Diophantische Approximationen und Hausdorff
Ma{\ss},} Math.~Sbornik, 36 (1929) 371--382.

\bibitem[KS]{KS} B.~Kra and J.~Schmeling, {\it Diophantine numbers, dimension
and Denjoy maps,} Acta Arithmetica (to appear).

\bibitem[L]{L} J.~Levesley, {\it A general inhomogeneous Jarnik-Besicovitch
theorem} J.~Number Theory 71 (1998) 65--80.

\bibitem[P]{Pesin} Ya.~Pesin, {\it Dimension theory in dynamical systems,}
University of Chicago Press, Chicago, 1997.

\bibitem[R]{R} Th.~Ruijgrok, {\it Periodic orbits in triangular billiards,} Acta
Physical Polonica B 22 (1991) 955--981.

\bibitem[S]{WS} W.  Schmidt, {\it Metrical theorems on fractional parts of
sequences}, TAMS 110 (1964), 493--518

\bibitem[T]{T} S.~Tabachnikov, {\it Billiards,} ``Panoramas et Syntheses'',
Soc.~Math.~France (1995).

\bibitem[Tr]{Tr} S.~Troubetzkoy, {\it Recurrence and periodic billiard orbits in
polygons,} preprint.

}

\end{thebibliography}
\end{document}